\documentclass{amsart}
\usepackage{eurosym}
\usepackage{amssymb,amsmath,amsfonts}
\usepackage{tikz}
\usepackage[all]{xy}
\usepackage[unicode=true]{hyperref}

\newtheorem{theorem}{Theorem}[section]
\newtheorem{corollary}[theorem]{Corollary}

\newtheorem{definition}[theorem]{Definition}
\newtheorem{example}[theorem]{Example}

\newtheorem{lemma}[theorem]{Lemma}
\newtheorem{remark}[theorem]{Remark}
\newtheorem{proposition}[theorem]{Proposition}
\def\Im{\textup{Im}}
\def\Ker{\textup{Ker}}
\def\Coker{\textup{Coker}}
\def\Tr{\textup{Tr}}

\newcommand{\cref}[1]{eq. (\ref{#1})}

\begin{document}
\title[Matrix invertible extensions]{Matrix invertible
extensions over commutative rings. Part I: general theory}
\author{Grigore C\u{a}lug\u{a}reanu, Horia F. Pop, Adrian Vasiu}

\date{Accepted for publication in final form in J.\ Pure Appl.\ Algebra on Nov. 18, 2024.}

\begin{abstract}
A unimodular $2\times 2$ matrix with entries in a commutative $R$ is called extendable (resp.\ simply extendable) if it extends to an invertible $3\times 3$ matrix (resp.\ invertible $3\times 3$ matrix whose $(3,3)$ entry is $0$). We obtain necessary and sufficient conditions for a unimodular $2\times 2$ matrix to be extendable (resp.\ simply extendable) and use them to study the class $E_2$ (resp.\ $SE_2$) of rings $R$ with the property that all unimodular $2\times 2$ matrices with entries in $R$ are extendable (resp.\ simply extendable). We also study the larger class $\Pi_2$ of rings $R$ with the property that all unimodular $2\times 2$ matrices of determinant $0$ and with entries in $R$ are (simply) extendable (e.g., rings with trivial Picard groups or pre-Schreier domains). Among Dedekind domains, polynomial rings over $\mathbb Z$ and Hermite rings, only the EDRs belong to the class $E_2$ or $SE_2$. If $R$ has stable range at most $2$ (e.g., $R$ is a Hermite ring or $\dim(R)\le 1$), then $R$ is an $E_2$ ring iff it is an $SE_2$ ring. 
\end{abstract}

\maketitle

\section{Introduction}

\label{S1}

Let $R$ be a commutative ring with identity; we denote by $U(R)$ its group of units, by $J(R)$ its Jacobson radical, and by $Pic(R)$ its Picard group. For $m,n\in\mathbb{N}=\{1,2,\ldots\}$, let $\mathbb{M}_{m\times n}(R)$ be the $R$-module of $m\times n$ matrices with entries in $R$; we view $\mathbb M_n(R):=\mathbb M_{n\times n}(R)$ as an $R$-algebra with identity $I_n$. Let ${GL}_{n}(R)$ be the general linear group of units of $\mathbb M_n(R)$, and let ${SL}_{n}(R):=\{M\in GL_n(R)|\det(M)=1\}$ be the special linear subgroup of ${GL}_{n}(R)$.

For a free $R$-module $F $, let $Um(F)$ be the set of 
\textsl{unimodular} elements of $F$, i.e., of elements $v\in F$ for which
there exists an $R$-linear map $L:F\rightarrow R$ such that $L(v)=1$; we
have an identity $U(R)=Um(R)$ of sets.

For $A\in \mathbb{M}_{n}(R)$, let $A^{T}$ be its transpose, let $\Tr(A)$ be its trace, let $\chi_{A}(\lambda )\in R[\lambda ]$ be its (monic) characteristic polynomial, and let $\nu_A:=\chi_A'(0)\in R$. E.g., if $n=3$, then $\chi_A(\lambda)=\lambda^3-\Tr(A)\lambda^2+\nu_A\lambda-\det(A)$. 

As we will be using pairs, triples and quadruples
extensively, the elements of $R^{n}$ will be denoted as $n$-tuples except
for the case of the $R$-linear map $L_A:R^{n}\rightarrow R^{n}$ defined 
by $A$, in which case they will be viewed as $n\times 1$ columns. 
Let $\Ker_{A}$, $\Coker_A$ and $\Im_{A}$ be the kernel, cokernel and image (respectively) of $L_{A}$. 

Recall that $B,C\in\mathbb M_{m\times n}(R)$ are said to be equivalent 
if there exist matrices $M\in {GL}_m(R)$ and $N\in GL_n(R)$ such that $C=MBN$; if $m=n$ and we can 
choose $N=M^{-1}$, then $B$ and $C$ are said to be similar. If all entries of $B-C$ belong to an ideal $I$ of $R$, then we say that $B$ and $C$ are congruent modulo $I$. By the reduction of $A\in\mathbb M_n(R)$ modulo $I$ we mean the image of $A$ in $\mathbb M_n(R/I)\cong\mathbb M_n(R)/\mathbb M_n(I)$.

We say that $B\in\mathbb M_{m\times n}(R)$ \textsl{admits diagonal reduction} if it is
equivalent to a matrix whose off diagonal entries are $0$ and whose diagonal
entries $b_{1,1},\ldots ,b_{s,s}$, with $s:=\min \{m,n\}$, are such that $b_{i,i}$ divides $b_{i+1,i+1}$ for all $i\in \{1,\ldots ,s-1\}$.

We will use the shorthand `iff' for `if and only if' in all that follows.

For $Q\in\mathbb M_3(R)$, let $\Theta(Q)\in \mathbb{M}_{2}(R)$ be obtained from $Q$ by removing the third row and the third column. If $Q\in GL_3(R)$, then $\Theta(Q)$ modulo each maximal ideal
of $R$ is nonzero, hence $\Theta(Q)\in Um\bigl(\mathbb{M}_2(R)\bigr)$. The rule $Q\rightarrow \Theta(Q)$ defines a
map 
\begin{equation*}
\Theta:=\Theta_R:{SL}_3(R)\rightarrow Um\bigl(\mathbb{M}_2(R)\bigr)
\end{equation*}and in this paper we study the image of $\Theta $ and in particular
the class of rings $R$ for which $\Theta$ is surjective (i.e., it
has a right inverse).

If $M=\left[ 
\begin{array}{cc}
a & b \\ 
c & d\end{array}\right]\in {GL}_{2}(R)$, let $\sigma(M):=\left[ 
\begin{array}{ccc}
a & b & 0\\ 
c & d & 0\\
0 & 0 & \det(M)^{-1}
\end{array}\right]\in {SL}_3(R)$. Clearly 
$\Theta\bigl(\sigma(M)\bigr)=M$, hence $GL_{2}(R)\subseteq {Im}(\Theta)$.

\begin{definition}\normalfont\label{def1}
We say that a matrix $A\in \mathbb{M}_2(R)$ is $SL_3$-extendable if there exists $A^{+}\in {SL}_3(R)$ such that $A=\Theta(A^+)$, and we call $A^{+}$ an $SL_3$-extension of $A$. If we can choose $A^{+}$ such that its $(3,3)$ is $0$, then we say that $A$ is simply $SL_3$-extendable and that $A^{+}$ is a simple $SL_3$-extension of $A$.
\end{definition}

As in this paper we do not consider $SL_n$-extensions with $n\ge 4$, in all that follows it will be understood
that ``extendable" means $SL_3$-extendable. Each extendable matrix is unimodular. Theorem \ref{TH4} contains three other equivalent characterizations of simply extendable matrices; e.g., a unimodular $2\times 2$ matrix admits diagonal reduction iff it is simply extendable.

The problem of deciding if $A\in \mathbb{M}_2(R)$ is simply extendable relates to classical studies of finitely generated stable free modules that aim to complete 
matrices in $\mathbb M_{n\times (n+m)}(R)$ whose $n\times n$ minors 
generate $R$ (e.g., see \cite{LL}), but it fits within the general problem of finding square matrices with prescribed entries and coefficients of characteristic polynomials. Concretely, if $A\in  Um\bigl(\mathbb{M}_2(R)\bigr)$ is simply extendable, its simple extensions $A^{+}$ have 5 prescribed entries (out of 9) and 2 prescribed coefficients (out of 3) of $\chi_{A^{+}}(\lambda )=\lambda ^{3}-\Tr(A)\lambda ^{2}+\nu_{A^{+}}\lambda -1$; the nonempty subsets 
\begin{equation*}
\nu(A):=\{\nu_{A^{+}}|A^{+}\; \textup{is a simple extension of}\; A\}\subseteq R
\end{equation*}are sampled in Examples \ref{EX5} and \ref{EX10}(1). Such general problems over fields have a long history 
and often complete results (cf.\ \cite{FL} and \cite{her}).

In general $\Theta$ is not surjective (see Theorem \ref{TH3}(4) and Example \ref{EX8}).

If $\det (A)=0$, then $A$ is extendable iff it is simply extendable. More generally, a matrix $A\in \mathbb{M}_2(R)$ is extendable iff its reduction modulo $R\det (A)$ is simply extendable (see Lemma \ref{L1}(1)). Based on this, Definition \ref{def1} leads to the study of 3 classes of rings to be named, in the
spirit of \cite{GH1} and \cite{lor}, using indexed letters.

\begin{definition}\normalfont\label{def2}
We say that $R$ is:

\medskip \textbf{(1)} a ${\Pi}_2$ ring, if each $A\in Um\bigl(\mathbb{M}_2(R)\bigr)$ with $\det(A)=0$ is extendable;

\smallskip \textbf{(2)} an ${E}_2$ ring, if each matrix in $ Um\bigl(\mathbb{M}_2(R)\bigr)$ is extendable (i.e., $\Theta$ is surjective);

\smallskip \textbf{(3)} an ${SE}_2$ ring, if each matrix in $ Um\bigl(\mathbb{M}_2(R)\bigr)$ is simply extendable.

\medskip
If moreover $R$ is an integral domain, we replace ring by a domain (so we speak about $\Pi_2$ domains, $E_2$ domains, etc.)
\end{definition}

We recall that $R$ is called an \textsl{elementary divisor ring}, abbreviated as \textsl{EDR}, if all matrices with entries in $R$ admit diagonal reduction. Equivalently, $R$ is
an \textsl{EDR} iff each matrix in $\mathbb{M}_2(R)$ is
equivalent to a diagonal matrix and iff each finitely presented $R $-module is a direct sum of cyclic $R$-modules (see \cite{LLS}, Cor.\ (3.7)
and Thm.\ (3.8); see also \cite{WW}, Thm.\ 2.1). 

Let $A\in Um\bigl(\mathbb{M}_2(R)\bigr)$. If $A$ admits diagonal reduction, let $M,N\in {GL}_2(R)$ be such that $MAN=\left[ 
\begin{array}{cc}
1 & 0\\ 
0 & \det(A)
\end{array}
\right]$; for $A^+:=\sigma(M^{-1})\left[ 
\begin{array}{ccc}
1 & 0 & 0\\ 
0 & \det(A) & 1\\
0 & -1 & 0
\end{array}\right] \sigma(N^{-1})$
one easily computes $\Theta(A^+)=M^{-1}MANN^{-1}=A$ (this is a particular case of Lemma \ref{L1}(2)), hence $A$ is simply extendable. We conclude that:

\begin{proposition}\label{PR1} 
Each \textsl{EDR} is an $SE_2$ ring.
\end{proposition}

The following characterizations of ${\Pi}_2$ rings are presented in Section \ref{S5}.

\begin{theorem}
\label{TH1} For a ring $R$ the following statements are equivalent:

\medskip \textbf{(1)} The ring $R$ is a ${\Pi}_{2}$ ring.

\smallskip \textbf{(2)} Each matrix in $Um\bigl(\mathbb{M}_2(R)\bigr)$ of zero
determinant is non-full, i.e., the product of two matrices of sizes $2\times
1$ and $1\times 2$ (equivalently, the $R$-linear map $L_A$ factors as a
composite $R$-linear map $R^2\rightarrow R\rightarrow R^2$).

\smallskip \textbf{(3)} For each matrix in $Um\bigl(\mathbb{M}_2(R)\bigr)$ of zero
determinant, one (hence both) of the $R$-modules $\Ker_A$ and $\Im_A$ is
isomorphic to $R$.

\smallskip \textbf{(4)} Each projective $R$-module of rank $1$ generated by
two elements is isomorphic to $R$.
\end{theorem}

For pre-Schreier domains see Section \ref{S2}. Recall that an integral domain $R$ is a pre-Schreier domain iff each matrix in $\mathbb{M}_{2}(R)$ of zero determinant
is non-full (see \cite{CP}, Thm.\ 1 or \cite{MR}, Lem.\ 1). Each pre-Schreier domain is a $\Pi _{2}$ domain by Theorem \ref{TH1}. Similarly, if $Pic(R)$ 
is trivial, then $R$ is a ${\Pi }_{2}$ ring.

The following notions introduced by Bass (see \cite{bas}), Shchedryk (see 
\cite{shc1}), and McGovern (see \cite{mcg1}) (respectively) will be often
used in what follows.

\begin{definition}\normalfont\label{def3}
Let $n\in\mathbb{N}$. Recall that $R$ has:

\medskip \noindent \textbf{(1)} \textsl{stable range} n and we write $sr(R)=n $ if $n$ is the smallest natural number with the property that each $(a_1,\ldots,a_n,b)\in Um(R^{n+1})$ is reducible, i.e., there exists $(r_1,\ldots,r_n)\in R^n$ such that $(a_1+br_1,\ldots,a_n+br_n)\in Um(R^n)$ 
(when there exists no $n\in\mathbb N$ such that $sr(R)=n$, then 
$sr(R):=\infty$ and the convention is $\infty>n$);

\smallskip \noindent \textbf{(2)} \textsl{(fractional) stable range 1.5} and 
we write $fsr(R)=1.5$ if for each $(a,b,c)\in Um(R^{3})$ with $c\neq 0$ 
there exists $r\in R$ such that $(a+br,c)\in Um(R^{2})$;

\smallskip\noindent \textbf{(3)} \textsl{almost stable range 1} and we
write $asr(R)=1$ if for each ideal $I$ of $R$ not contained in $J(R)$, $sr(R/I)=1$.
\end{definition}

Stable range type of conditions on all suitable unimodular tuples 
of a ring date back at least to Kaplansky. For instance, the essence of 
\cite{kap}, Thm.\ 5.2 can be formulated in the language of this paper 
as follows: each triangular matrix in $Um\bigl(\mathbb{M}_2(R)\bigr)$ is
simply extendable iff for each $(a,b,c)\in Um(R^{3})$, there exists $(e,f)\in R^2$ such that 
$(ae,be+cf)\in Um(R^2)$. The general (nontriangular) form of this reformulation is: a ring $R$ is an $SE_2$ 
ring iff for each $(a,b,c,d)\in Um(R^4)$, there exists $(e,f)\in R^2$ such that $(ae+cf,be+df)\in Um(R^2)$ (see Theorem \ref{TH4}).

We recall that if $R$ has (Krull) dimension $d$, then $sr(R)\le d+1$ (see \cite{bas}, Thm.\ 11.1 for the noetherian case and see \cite{hei}, Cor.\ 2.3 for the general case). Also, if $R$ is a finitely generated algebra over a finite field of dimension $2$ 
or if $R$ is a polynomial algebra in 2 indeterminates over a field that is algebraic 
over a finite field, then $sr(R)\le 2$ (see \cite{VS}, 
Cors.\ 17.3 and 17.4).
 
Each ${SE}_{2}$ ring is an ${E}_2$ ring, but we do not know when
the converse is true. However, we show that there exist $2\times 2$ matrices that are extendable but are not simply extendable (see Example \ref{EX9}). Moreover, we have (see Section \ref{S6}): 

\begin{theorem}\label{TH2} 
If $sr(R)\le 2$, then the extendable and simply extendable properties on a matrix in $\mathbb M_2(R)$ are equivalent (hence $R$ is an $SE_2$ ring iff it is an $E_2$ ring).
\end{theorem}

Based on Theorem \ref{TH1}, in Section \ref{S6} we prove the following theorem.

\begin{theorem}
\label{TH3} Let $R$ be an integral domain of dimension 1. Then the following properties hold:

\medskip \textbf{(1)} Each matrix in $Um\bigl(\mathbb{M}_2(R)\bigr)$ with nonzero determinant is
simply extendable.

\smallskip \textbf{(2)} Each triangular matrix in $Um\bigl(\mathbb{M}_2(R)\bigr)$ is
simply extendable.

\smallskip \textbf{(3)} The ring $R$ is a ${\Pi}_2$ domain iff it is an ${SE}_2$ (or an $E_2$) domain and iff $Pic(R)$ is trivial. 

\smallskip \textbf{(4)} Assume $R$ is a Dedekind domain. The ring $R$ is a ${\Pi}_2$ domain iff it is a principal ideal domain (\textsl{PID}).
\end{theorem}

Recall that $R$ is a \textsl{Hermite} ring in the sense of Kaplansky, if $RUm(R^{2})=R^{2}$, equivalently
if each $1\times 2$ matrix with entries in $R$ admits diagonal reduction. If $R$ is a Hermite ring, then $sr(R)\in \{1,2\}$ (see \cite{MM}, Prop.\ 8(i); see also \cite{zab1}, Thm.\ 2.1.2) and a simple induction on $n\in\mathbb N$ gives that $RUm(R^{n})=R^{n}$; thus $RUm\bigl(\mathbb{M}_2(R)\bigr)=\mathbb M_2(R)$. It follows that a Hermite ring $R$ is an EDR iff each matrix in $Um\bigl(\mathbb{M}_2(R)\bigr)$ admits diagonal reduction (equivalently, it is simply completable, see Theorem \ref{TH4}). From the last two sentences and Theorem \ref{TH2} we conclude:

\begin{corollary}\label{C1}
Let $R$ be a Hermite ring. Then $R$ is an EDR iff it is an $E_2$ ring and iff it is an $SE_2$ ring.
\end{corollary}

For almost stable range 1 we have the following applications (see Section \ref{S6}):

\begin{corollary}
\label{C2} Assume that $asr(R)=1$. Then the following properties hold:

\medskip \textbf{(1)} Each triangular matrix in $Um\bigl(\mathbb{M}_2(R)\bigr)$ is
simply extendable.

\smallskip \textbf{(2) (McGovern)} If $R$ is a Hermite ring, then $R$ is an 
\textsl{EDR}.
\end{corollary}

Corollary \ref{C2}(2) was first obtained in \cite{mcg1},
Thm.\ 3.7. A second proof of Corollary \ref{C2}(2) is presented in  Remark \ref{rem5}.

Part II studies determinant liftable $2\times 2$ matrices that generalize simply extendable matrices\footnote{A matrix $A\in Um\bigl(\mathbb{M}_2(R)\bigr)$ will be called determinant liftable if there exists $B\in Um\bigl(\mathbb{M}_2(R)\bigr)$ congruent to $A$ modulo $R\det(A)$ and $\det(B)$=0.} and proves that each $J_{2,1}$ domain introduced in \cite{lor} is an EDR. Part III has applications to B\'{e}zout rings (i.e., to rings whose finitely generated ideals are principal). Part IV contains universal and stability properties, complements and open problems. Parts I to IV split the manuscript \url{https://arxiv.org/abs/2303.08413}.

\section{Basic terminology and properties}\label{S2}

In what follows we will use without extra comments the following two basic properties. For $(a,b,c)\in R^3$ we have $(a,bc)\in Um(R^2)$ iff 
$(a,b),(a,c)\in Um(R^2)$. If $(a,b)\in Um(R^2)$ and $c\in R$, then $a$ divides $bc$ iff $a$ divides $c$.

A ring $R$ is called \textsl{pre-Schreier}, if every nonzero element $a\in R$
is \textsl{primal}, i.e., if $a$ divides a product $bc$ of elements of $R$,
there exists $(d,e)\in R^2$ such that $a=de$, $d$ divides $b$ and $e$
divides $c$. Pre-Schreier domains were introduced by Zafrullah in \cite{zaf}.
A pre-Schreier integrally closed domain was called a \textsl{Schreier}
domain by Cohn in \cite{coh}. Every \textsl{GCD} domain (in particular,
every B\'{e}zout domain) is Schreier (see \cite{coh}, Thm.
2.4). Products of pre-Schreier domains and quotients of \textsl{PID}s are pre-Schreier rings.

In an integral domain, an irreducible element is primal iff it is
a prime. Thus an integral domain that has irreducible elements which are not
prime, such as each noetherian domain which is not a UFD, is not
pre-Schreier.

The \textsl{inner rank} of an $m\times n$ matrix over a ring is defined as
the least positive integer $r$ such that it can be expressed as the product
of an $m\times r$ matrix and an $r\times n$ matrix; over fields, 
this notion coincides with the usual notion of rank. A square matrix
is called \textsl{full} if its inner rank equals its order, and \textsl{non-full} otherwise. A $2\times 2$ matrix is non-full iff its
inner rank is $1$, (i.e., it has a column-row decomposition).

We consider the subsets of $Um(R^3)$: 
\begin{equation*}
T_{3}(R):=Um (R^{2})\times (R\setminus \{0\})\;\text{\textup{and}}\;J_{3}(R):=Um (R^{2})\times \bigl(R\setminus J(R)\bigr).
\end{equation*}

\begin{proposition}
\label{PR2} We have $fsr(R)=1.5$ iff for each $(a,b,c)\in T_{3}(R)$ there exists $r\in R$ such that $(a+br,c)\in Um(R^{2})$.
\end{proposition}

\begin{proof}
The `only if' part is clear. To check the `if' part, let $(a,b,c)\in
Um(R^{3})$ with $c\neq 0$ and let $(x,y,z)\in R^3$ be such that $ax+by+cz=1$.
Thus $(a,by+cz,c)\in T_{3}(R)$ and hence there exists $r\in R$ such that $(a+byr+czr,c)\in Um(R^{2})$. This implies $(a+byr,c)\in Um(R^{2})$, thus $fsr(R)=1.5$.
\end{proof}

\begin{proposition}
\label{PR3} We have $asr(R)=1$ iff for each $(a,b,c)\in J_{3}(R)$ there exists $r\in R$ such that $(a+br,c)\in Um(R^{2})$.
\end{proposition}

\begin{proof}
See \cite{mcg1}, Thm.\ 3.6 for the `only if' part. For the `if' part, for $I$ an ideal of $R$ not contained in $J(R)$ we check that $sr(R/I)=1$. If $(a,b)\in R^2$ is such that $(a+I,b+I)\in Um\bigl((R/I)^{2}\bigr)$, let
$(d,e)\in R^2$ and $c\in I$ be such that $ad+be+c=1$. If $c\notin J(R)$,
then for $(f,g):=(c,c)\in \bigl(I\setminus J(R)\bigr)\times I$ we have $(a,be+g,f)\in J_{3}(R)$. If $c\in J(R)$, then $ad+be=1-c\in U(R)$, hence $(a,be)\in Um(R^2)$ and for $(f,g)\in\bigl(I\setminus J(R)\bigr)\times \{0\}$ 
we have $(a,be+g,f)\in J_{3}(R)$. If $r\in R$ is such that $\bigl(a+(be+g)r,f\bigr)\in Um(R^{2})$, then $a+I+(b+I)(er+I)\in U(R/I)$; so $sr(R/I)=1$.
\end{proof}

\begin{corollary}
\label{C3} \textbf{(1)} If $sr(R)=1$, then $fsr(R)=1.5$.

\smallskip \textbf{(2)} If $fsr(R)=1.5$, then $asr(R)=1$.

\smallskip \textbf{(3)} If $asr(R)=1$, then $sr(R)\le 2$.
\end{corollary}

\begin{proof}
If $sr(R)=1$, then for each $(a,b,c)\in T_3(R)$, there exists $r\in R$ such that $a+rb\in U(R)$, so $(a+rb,c)\in Um(R^2)$, hence $fsr(R)=1.5$ by Proposition \ref{PR2}. So part (1) holds. As $J_{3}(R)\subseteq T_{3}(R)$, part (2) follows from Propositions \ref{PR2} and \ref{PR3}.

To check part (3), it suffices to show that each $(a,b,c)\in Um(R^3)$ is
reducible. If $b\notin J(R)$, then $R/Rb$ has stable range $1$ and $(a+Rb,c+Rb)\in Um\bigl((R/Rb)^2\bigr)$; hence there exists $r\in R$ such that $a+cr+Rb\in U(R/Rb)$ and thus for $(r_1,r_2):=(r,0)$ we have $(a+cr_1,b+cr_2)\in Um(R^2)$. If $b\in J(R)$, then $(a,b,c)\in Um(R^3)$
implies that $(a,c),(a,b+c)\in Um(R^2)$ and thus for $(r_1,r_2):=(0,1)$ we
have $(a+cr_1,b+cr_2)\in Um(R^2)$. We conclude that $(a,b,c)$ is reducible.
\end{proof}

Corollary \ref{C3}(3) was first obtained in \cite{mcg1}, Thm.\ 3.6.

We have the following `classical' \textsl{units and unimodular interpretations}.

\begin{proposition}
\label{PR4} \textbf{(1)} For $n\in\mathbb{N}$, we have $sr(R)\le n$ iff for each $b\in R$ the
reduction modulo $Rb$ map of sets $Um(R^n)\rightarrow Um\bigl((R/Rb)^n\bigr)$ is
surjective.

\smallskip \textbf{(2)} We have $fsr(R)=1$ iff for
all $(b,c)\in R^2$ with $c\neq 0$, the homomorphism $U(R/Rc)\rightarrow
U\bigl(R/(Rb+Rc)\bigr)$ is surjective.

\smallskip \textbf{(3)} We have $asr(R)=1$ iff for each $(b,c)\in R^2$
with $c\notin J(R)$, the homomorphism $U(R/Rc)\rightarrow U\bigl(R/(Rb+Rc)\bigr)$ is
surjective.
\end{proposition}

\begin{proof}
Parts (1) and (2) follow from definitions. To check the `if' part of (3), let $(a,b,c)\in J_{3}(R)$. Then $a+Rb+Rc$
is a unit of $R/(Rb+Rc)$ and thus is the image of a unit of $R/Rc$, which is
of the form $a+br+Rc$ with $r\in R$. Hence $(a+br,c)\in Um(R^{2})$. From
this and Proposition \ref{PR3} it follows that $asr(R)=1$. To check the `only if' part of (3), let $a+Rb+Rc\in U\bigl(R/(Rb+Rc)\bigr)$. Thus $(a+Rc,b+Rc)\in Um\bigl((R/Rc)^2\bigr)$. As $Rc\not\subseteq J(R)$ and we assume $asr(R)=1$, it follows that $sr(R/Rc)=1$, hence there exists $r\in R$ such that 
$a+br+Rc\in U(R/Rc)$ maps to $a+Rb+Rc\in U\bigl(R/(Rb+Rc)\bigr)$. Thus the homomorphism $U(R/Rc)\rightarrow U\bigl(R/(Rb+Rc)\bigr)$ is surjective.
\end{proof}

\begin{corollary}
\label{C5} Assume that $sr(R)\le 4$. If $R$ is an $E_{2}$ (resp.\ $SE_2$) ring, then $R/Ra$ is an $E_{2}$ (resp.\ $SE_2$) ring for all $a\in R$.
\end{corollary}

\begin{proof}
Let $\bar A\in Um\bigl(\mathbb M_2(R/Ra)\bigr)$. As $sr(R)\le 4$, there exists $A\in Um\bigl(\mathbb M_2(R)\bigr)$ whose reduction modulo $Ra$ is $\bar A$ by Proposition \ref{PR4}(1) applied to $n=4$. Let $A^+$ be an extension (resp.\ a simple extension) of $A$; its reduction modulo $Ra$ is an extension (resp.\ a simple extension) of $\bar A$. So $R/Ra$ is an $E_2$ (resp.\ $SE_2$) ring.\end{proof}

\begin{example}
\normalfont
\label{EX1} Let $(a,b,c)\in Um(R^3)$ with $(b,c)\in Um(R^2)$. Writing $a-1=be+cf$ with $e,f\in R$, for $r:=-e$ we have $(a+br,c)\in Um(R^2)$.
\end{example}

\begin{example}
\normalfont\label{EX2} Suppose $R$ is a semilocal ring. For each $b\in R$, the homomorphism $U(R)\rightarrow U(R/Rb)$ is
surjective, hence $sr(R)=1$ by Proposition \ref{PR4}(1).
\end{example}

\begin{example}
\normalfont\label{EX3} Let $R$ be a noetherian domain of dimension 1. For $(b,c)\in R^2$ with $c\neq 0$, the rings $R/Rc$ and $R/(Rb+Rc)$ are artinian, hence the homomorphism $U(R/Rc)\rightarrow U\bigl(R/(Rb+Rc)\bigr)$ is surjective. From
this and Proposition \ref{PR4}(2), it follows that $fsr(R)=1.5$;
thus $sr(R)\le 2$ (see Corollary \ref{C3}(2) and (3)).
\end{example}

An argument similar to the one of Example \ref{EX3} shows that each B\'{e}zout domain which is a filtered union of Dedekind domains has stable range
1.5.

\begin{example}
\normalfont\label{EX4} We have $sr(\mathbb{Z})=2$ and $fsr(\mathbb Z)=1.5$. The ring $\mathbb{Z}[x]/(x^{2})$ has almost stable range 1 but does not have stable range 1.5. For a finite field $F$ and indeterminaes $X,Y$, $sr(F[X,Y])=sr(F[X])=2$; so $F[X,Y]$ does not have almost stable range 1. Hence, each possible converse of Corollary \ref{C3} does not hold.
\end{example}

\section{Projective modules}\label{S3}

For a projective $R$-module $P$ of rank $1$, let $[P]\in {Pic}(R)$ be its class in the
Picard group. Let ${Pic}_{2}(R)$ be the subgroup of ${Pic}(R)
$ generated by classes $[P]$ with $P$ generated by $2$ elements. Let ${Pic}_{2}(R)[2]$ be the subgroup of ${Pic}_{2}(R)$ generated
by classes $[P]$ with $P\oplus P\cong R^2$ (so $2[P]=[R]$ and $P$ is generated by $2$ elements).

For $(\alpha _{1},\ldots ,\alpha _{n})\in R^{n}$, let 
${Diag}(\alpha_{1},\ldots ,\alpha _{n})\in \mathbb{M}_{n}(R)$ be the diagonal matrix 
whose $(i,i)$ entry is $\alpha _{i}$ for all $i\in \{1,\ldots ,n\}$.

If $n\geq 2$ we do not have unit interpretations of the stable range $n$ similar to Proposition \ref{PR4} but this is replaced by standard projective modules considerations recalled here in the form required in the sequel. If $P_1$ and $P_2$ are two projective $R$-modules of rank $1$ such that $P_1\oplus P_2\cong R^2$, then by taking determinants it follows that $P_1\otimes_R P_2\cong R$, hence $[P_1]=-[P_2]\in Pic(R)$; thus $P_1\cong R$ iff $P_2\cong R$. 

\begin{lemma}\label{L0}
Let $A\in Um\bigl(\mathbb{M}_2(R)\bigr)$. Let $\bar{R}:=R/R\det(A)$ and let $\bar A$ be the reduction of $A$ modulo $R\det(A)$. Then $\Ker_{\bar{A}}$, $\Im_{\bar{A}}$ and $\Coker_{\bar{A}}$ are projective $\bar{R}$-modules of rank $1$ generated by two elements and we have $[\Ker_{\bar{A}}]=[\Coker_{\bar{A}}]=-[\Im_{\bar{A}}]\in {Pic}_2(R)$ (thus $\Ker_{\bar{A}}\cong\bar{R}$ iff $\Im_{\bar{A}}\cong\bar{R}$ and iff $\Coker_{\bar{A}}\cong\bar{R}$).
\end{lemma}

\begin{proof}
As $\Coker_A$ is annihilated by $\det(A)$, we can view it as an $\bar{R}$-module isomorphic to $\Coker_{\bar A}$. Locally in the Zariski
topology of the spectrum of $\bar{R}$, one of the entries of $\bar{A}$ is a unit and
hence the matrices $\bar{A}$ and ${Diag}(1,0)$ are equivalent. Thus $\Im_{\bar{A}}$ and $\Coker_{\bar{A}}$ are projective $\bar{R}$-modules of rank $1$ generated by two elements and we have two (split) short exact sequences of projective $\bar{R}$-modules $0\rightarrow\Im_{\bar{A}}\rightarrow \bar{R}^2\rightarrow \Coker_{\bar{A}}\rightarrow 0$ and $0\rightarrow\Ker_{\bar{A}}\rightarrow
\bar{R}^{2}\rightarrow\Im_{\bar{A}}\rightarrow 0$. From the existence of $\bar{R}$-linear isomorphisms $\bar{R}^{2}\cong\Ker_{\bar{A}}\oplus\Im_{\bar{A}}\cong\Ker_{\bar{A}}\oplus \Coker_{\bar{A}}$ it follows that $\Ker_{\bar{A}}$ and $\Coker_{\bar{A}}$ are isomorphic to the dual of $\Im_{\bar{A}}$ and the lemma follows.
\end{proof}

Nest we exemplify the connection between projective $R$-modules and reducibility of $n+1$-tuples with entries in $R$.

\begin{example}\normalfont\label{EX4.1}
For $n\in\mathbb{N}$, let $(a_1,\ldots,a_n,b)\in Um(R^{n+1})$ with 
$b$ a nonzero divisor. We consider short exact sequences of $R$-modules $0\rightarrow R\xrightarrow{b} R\xrightarrow{\pi} R/Rb\rightarrow 0$ and $0\rightarrow
Q\rightarrow R^n\xrightarrow{f} R/Rb\rightarrow 0$ where $\pi$ is the natural quotient map and $f$ maps the
elements of the standard basis of $R^n$ to $a_1+Rb,\ldots,a_n+Rb$ and $Q:=Ker(f)$. If $Q^+\rightarrow R^n$ and $Q^+\rightarrow R$ define the pullback of $f$ and $\pi$, then we have short exact sequences $0\rightarrow Q\rightarrow
Q^+\rightarrow R$ and $0\rightarrow R\rightarrow Q^+\rightarrow
R^n\rightarrow 0$ which imply that $Q^+$ is a free $R$-module of rank $n+1$
and $Q$ is a projective $R$-module of rank $n$ generated by $n+1$ elements; moreover, if $n=1$, then $Q\cong R$.

We consider $R$-linear maps $g:R^{n}\rightarrow R$ such that $\pi \circ g=f$. Then $(a_{1},\ldots ,a_{n},b)$ is reducible iff we can choose $g $ to be surjective. Thus, if $sr(R)\le n$, then we can choose $g$ to be surjective.

If $A\in \mathbb{M}_{n}(R)$ is equivalent to ${Diag}(1,1,\ldots,1,b) $ and $\Im_A$ is the $R$-submodule $Q$ of $R^{n}$ (hence $Q\cong R^n$), then we can choose $g$ to be surjective and hence $(a_{1},\ldots ,a_{n},b)$ is reducible.

In this paragraph we assume that $(a_{1},\ldots ,a_{n},b)$ is reducible and
that $g$ is chosen to be surjective. Then ${Ker}(g)$ is a projective $R$-module of rank $n-1$ generated by $n$ elements, we have a short exact
sequence $0\rightarrow {Ker}(g)\rightarrow Q\rightarrow Rb\rightarrow 0$,
and $0\rightarrow Q\rightarrow R^{n}\xrightarrow{f}R/Rb\rightarrow 0$
is the direct sum of the two projective resolutions $0\rightarrow R\xrightarrow{b}R\xrightarrow{\pi}R/Rb\rightarrow 0$ and $0\rightarrow {Ker}(g)\rightarrow {Ker}(g)\rightarrow 0\rightarrow 0$. If $n=2$, then ${Ker}(g)\cong R$.
Thus if $n=2$ and the $R$-submodule $Q$ of $R^{2}$ is $\Im_{A}$ for some $A\in 
\mathbb{M}_{2}(R)$, then $A$ is equivalent to ${Diag}(1,b)$.
Similarly, if $n\geq 3$ and ${Ker}(g)\cong R^{n-1}$, then $Q$ is a free $R$-module of rank $n$, so if the $R$-submodule $Q$ of $R^{n}$ is $\Im_{A}$
for some $A\in \mathbb{M}_{n}(R)$, then $A$ is equivalent to ${Diag}(1,1,\ldots,1,b)$.\end{example}

\section{Criteria on extending $2\times 2$ matrices} \label{S4}

\begin{lemma}\label{L1}
The following hold:

\medskip \textbf{(1)} A matrix $A\in \mathbb{M}_2(R)$ is extendable iff its reduction modulo $R\det (A)$ is simply extendable. Thus, if $\det (A)=0$, then $A$ is extendable iff it is simply extendable.

\smallskip \textbf{(2)} For $M,N\in {GL}_2(R)$ and $Q\in {SL}_3(R)$ we have an identity 
\begin{equation*}
\Theta\bigl(\sigma(M)Q\sigma(N)\bigr)=M\Theta(Q)N
\end{equation*}and the $(3,3)$ entries of $Q$ and $\sigma(M)Q\sigma(N)$ 
are equivalent (thus one such entry is $0$ iff
the other entry is $0$). Also, $\Theta(Q^{T})=\Theta(Q)^{T}$.

\smallskip \textbf{(3)} The fact that a matrix $A\in \mathbb{M}_2(R)$ is
(simply) extendable depends only on the equivalence class $[A]\in {GL}_2(R)\backslash \mathbb{M}_2(R)/GL_2(R)$. Moreover, $A$ is (simply) extendable iff so is $A^{T}$. So ${Im}(\Theta)$ is stable under transposition and
equivalence.
\end{lemma}

\begin{proof}
To check the `only if' part of part (1), let $A\in\mathbb M_2(R)$ be extendable, with $A^+\in SL_3(R)$ an extension of it. If $A^+_0$ is obtained from $A^+$ by replacing the $(3,3)$ entry with $0$, then the reductions of $A^+$ and $A^+_0$ modulo $R\det(A)$ have the same determinant $1$, and it follows that $A$ modulo $R\det(A)$ is simply extendable. To check the `if' part of (1), let $B\in\mathbb M_3(R)$ be such that $\Theta(B)=A$ and its reduction modulo $R\det(A)$ is a simple extension of the reduction of $A$ modulo $R\det(A)$. Let $w\in R$ be such that $\det(B)=1+w\det(A)$. If $A^+\in\mathbb M_3(R)$ is obtained from $B$ by subtracting $w$ from its $(3,3)$ entry, then $A^+$ is an extension of $A$ as $\det(A^+)=\det(B)-w\det(A)=1$. Thus part (1) holds. Part (2) is a simple computation, while part (3) follows directly from part (2).\end{proof}

\begin{corollary}\label{C6} We consider the following two statements on $R$.

\medskip
{\bf (1)} For each $a\in R$, $R/Ra$ is a $\Pi_2$ ring.

\smallskip
{\bf (2)} The ring $R $ is an $E_2$ ring.

\medskip
Then $(1)\Rightarrow (2)$. If $sr(R)\le 4$, then $1\Leftrightarrow (2)$.
\end{corollary}

\begin{proof}
To prove that $(1)\Rightarrow (2)$, let $A\in Um\bigl(\mathbb M_2(R)\bigr)$. Its reduction modulo $R\det(A)$ has zero determinant and hence it is simply extendable as $R/R\det(A)$ is a $\Pi_2$ ring. From Lemma \ref{L1}(1) it follows that $A$ is extendable. Thus $R$ is an $E_2$ ring. Hence $(1)\Leftrightarrow (2)$. If $sr(R)\le 4$, then $(2)\Rightarrow (1)$ by Corollary \ref{C5}, and hence $(1)\Leftrightarrow (2)$.\end{proof}

\begin{theorem}\label{TH4}
For $A=\left[ 
\begin{array}{cc}
a & b \\ 
c & d\end{array}\right]\in Um\bigl(\mathbb M_2(R)\bigr)$ the following statements are equivalent:

\medskip \textbf{(1)} The matrix $A$ is equivalent to the diagonal matrix ${Diag}\bigl(1,\det(A)\bigr)$.

\smallskip \textbf{(2)} The matrix $A$ is simply extendable.

\smallskip \textbf{(3)} There exists $(e,f)\in R^2$ such that $(ae+cf,be+df)\in Um(R^{2})$ (note that $(e,f)\in Um(R^2)$).

\smallskip \textbf{(4)} There exists $(x,y,z,w)\in R^4$ such that $ax+by+cz+dw=1$ and the matrix $\left[ 
\begin{array}{cc}
x & y \\ 
z & w\end{array}\right] $ is non-full.
\end{theorem}

\begin{proof}
For $(e,f,s,t)\in R^4$ we have two identities for the determinant
\begin{equation}  \label{EQ1}
\det \left[ 
\begin{array}{ccc}
a & b & f \\ 
c & d & -e \\ 
-t & s & 0\end{array}\right] =(be+df)t+(ae+cf)s=a(es)+b(et)+c(fs)+d(ft).
\end{equation}
The equivalence $(2)\Leftrightarrow (3)$ follows from the first identity. The implication $(1)\Rightarrow (2)$ was checked before Proposition \ref{PR1}. 

To show that $(3)\Rightarrow (1)\wedge (4)$, let $(e,f)\in
R^2$ be such that there exists $(s,t)\in Um(R^{2})$ with $(ae+cf)s+(be+df)t=1$. Then $(e,f),(s,t)\in Um(R^{2})$ and thus there exists $M\in {SL}_{2}(R)$ whose first row is $[e\;f]$ and there exists $N\in {SL}_{2}(R)$ whose first column is $[s\;t]^{T}$. The matrix $MAN$ has
determinant $\det (A)$ and its $(1,1)$ entry is $1$, thus it is equivalent to
the matrix ${Diag}\bigl(1,\det (A)\bigr)$ and therefore statement $(1)$ holds by the transitivity of the equivalent relation. For $(x,y,z,w):=(es,et,fs,ft)$, then $\left[ 
\begin{array}{cc}
x & y \\ 
z & w\end{array}\right]=\left[ 
\begin{array}{c}
e\\ 
f\end{array}\right]\left[ 
\begin{array}{cc}
s & t\end{array}\right]$ is non-full and from the second identity we get that $ax+by+cz+dw=1$, thus statement (4) holds. If statement (4) holds, let $(e,f,s,t)\in R^{4}$ be such that $(x,y,z,w)=(es,et,fs,ft)$; as $ax+by+cz+dw=1$, the determinant is $1$ by the second identity, hence $A$ is simply extendable.
\end{proof}

\begin{remark}\label{rem1}
\normalfont
Referring to Theorem \ref{TH4}, as $(2)\Leftrightarrow (3)$ and as a $2\times 2$ matrix has a
(simple) extension iff its transpose has it, it follows that
$A$ is simply extendable iff there exists $(e^{\prime
},f^{\prime })\in R^2$ such that $(ae^{\prime }+bf^{\prime },ce^{\prime
}+df^{\prime })\in Um(R^{2})$.
\end{remark}

\begin{example}
\normalfont\label{EX5} 
If $A^+=\left[ 
\begin{array}{ccc}
a & b & f \\ 
c & d & -e \\ 
-t & s & 0\end{array}\right]$ is a simple extension of $A=\left[ 
\begin{array}{cc}
a & b \\ 
c & d\end{array}\right]\in Um\bigl(\mathbb{M}_2(R)\bigr)$, then the characteristic polynomial $\chi_{A^+}$ of $A^+$ is of the form 
\begin{equation*}
x^3-{Tr}(A)x^2+\nu_{A^+}x-1
\end{equation*}
(see Section \ref{S1} for ${Tr}(A)=a+d$ and $\nu_{A^+}=\det(A)+es+ft$).
Thus the set of characteristic polynomials of simple extensions of $A$ is
in bijection to the subset $\nu(A)\subseteq R$ introduced in Section \ref{S1}
and we have 
\begin{equation*}
\nu(A)=\{\det(A)+es+ft|(e,f,s,t)\in R^4,\;a(es)+b(et)+c(fs)+d(ft)=1\}.
\end{equation*}

If $d\in R$, then $\nu\bigl(Diag(1,d)\bigr)=\{d+es+ft|(e,f,s,t)\in R^4,\ es+dft=1\}$ is equal to $\{d+1-(d-1)ft|(f,t)\in R^2\}=2+R(d-1)$.

Concretely, if $R=\mathbb{Z}$, $\{(8+11k,-5-7k)|k\in\mathbb{Z}\}$ is the solution set of the equation $7es+11ft=1$ in 
$es$ and $ft$, thus $\nu\bigl({Diag}(7,11)\bigr)=\{80+4k|k\in\mathbb Z\}=4\mathbb{Z}$.\end{example}

\begin{corollary}
\label{C7} \textbf{(1)} The ring $R$ is an ${SE}_2$ ring iff for
each $(a,b,c,d)\in Um(R^{4})$ statement (3) (or (4)) of Theorem \ref{TH4} holds.

\smallskip \textbf{(2)} Each semilocal ring is an ${SE}_2$ ring.
\end{corollary}

\begin{proof}
Part (1) follows from definitions and Theorem \ref{TH4}. For part (2), based on part (1) it suffices to show that for each $(a,b,c,d)\in Um(R^{4})$ there exists $(e,f)\in R^2$ such that $(ae+cf,be+df)\in Um(R^2)$. To prove this we can replace $R$ by $R/J(R)$; thus $R$ is a finite product of fields. By considering factors of $R$ that are fields, we can assume that $R$ is a field, in which case the existence of $(e,f)$ is clear.
\end{proof}

\begin{corollary}
\label{C8} For $(a,b,c,d)\in Um(R^{4})$ the following properties hold:

\medskip \textbf{(1)} The matrix $\left[ 
\begin{array}{cc}
a & b \\ 
c & d\end{array}\right] $ is extendable iff there exists $(e,f)\in R^2$ such that $(ae+cf,be+df,ad-bc)\in Um(R^3)$.

\smallskip \textbf{(2)} The matrix $\left[ 
\begin{array}{cc}
a & b \\ 
c & d\end{array}\right] $ is simply extendable iff there exists $(e,f)\in Um(R^2)$ such that $(ae+cf,be+df,ad-bc)\in Um(R^3)$
\end{corollary}

\begin{proof}
Part (1) follows from Lemma \ref{L1}(1) and Theorem \ref{TH4}. The `only if' of part (2) follows from Theorem \ref{TH4}. We are left to prove that if $(e,f)\in Um(R^2)$ is such that $(ae+cf,be+df,ad-bc)\in Um(R^3)$, then $A$ is simply extendable. Based on Theorem \ref{TH4} it suffices to show
that in fact we have $(ae+cf,be+df)\in Um(R^{2})$. Let $I:=R(ae+cf)+R(be+df)$ and $\mathfrak{m}$ a maximal ideal of $R$. If $ad-bc\in \mathfrak{m}$, then $(ae+cf,be+df,ad-bc)\in Um(R^{3})$ implies that 
$I\not\subseteq \mathfrak{m}$. If $ad-bc\notin \mathfrak{m}$, then $A$ modulo $\mathfrak{m}$ is invertible, hence $I\subseteq \mathfrak{m}$ iff $Re+Rf\subseteq \mathfrak{m}$; from this and $(e,f)\in Um(R^{2})$ we infer that 
$I\not\subseteq \mathfrak{m}$. As $I$ is not contained in any maximal ideal of 
$R$, it follows that $(ae+cf,be+df)\in Um(R^{2})$.
\end{proof}

\begin{corollary}
\label{C9} The following properties hold:

\medskip \textbf{(1)} The ring $R$ is an ${E}_2$ ring iff for
each $(a,b,c,d)\in Um(R^{4})$ there exists $(e,f)\in R^2$ such that $(ae+cf,be+df,ad-bc)\in Um(R^3)$.

\smallskip \textbf{(2)} The ring $R$ is an ${SE}_2$ ring iff for
each $(a,b,c,d)\in Um(R^{4})$ there exists $(e,f)\in Um(R^2)$ such that $(ae+cf,be+df,ad-bc)\in Um(R^3)$.
\end{corollary}

\begin{proof}
Both parts follow directly from definitions and the corresponding two parts of Corollary \ref{C8}.
\end{proof}

\begin{example}
\normalfont\label{EX6} Let $A=\left[ 
\begin{array}{cc}
a & b \\ 
c & d\end{array}\right]\in Um\bigl(\mathbb{M}_2(R)\bigr)$. In many simple cases one can easily
prescribe $(e,f)\in R^{2}$ such that $(ae+cf,be+df)\in Um(R^2)$ and hence conclude that $A$ is simply extendable. We include four such cases as follows.

\textbf{(1)} If $\{a,b,c,d\}\cap U(R)\neq \emptyset $, then we can take $(e,f)\in R^2$
such that $\{e,f\}=\{0,1\}$. E.g., if $a\in U(R)$, then $\left[ 
\begin{array}{ccc}
a & b & 0 \\ 
c & d & -1 \\ 
0 & a^{-1} & 0\end{array}\right] $ is a simple extension of $A$.

\textbf{(2)} If $\{(a,b),(a,c),(b,d),(c,d)\}\cap Um(R^{2})\neq \emptyset $,
then we can take $(e,f,e^{\prime },f^{\prime })\in R^4$ such that $1\in
\{ae+cf,be+df,ae^{\prime }+bf^{\prime },ce^{\prime }+df^{\prime }\}$. E.g., 
if $(a,b)\in Um(R^{2})$ and $s,t\in R$ are such that $as+bt=1$,
then $\left[ 
\begin{array}{ccc}
a & b & 0 \\ 
c & d & -1 \\ 
-t & s & 0\end{array}\right] $ is a simple extension of $A$.

\textbf{(3)} If at least two of the entries $a$, $b$, $c$ and $d$ are in $J(R)$ (e.g., they are $0$), then either $\{(a,b),(a,c),(b,d),(c,d)\}\cap Um(R^{2})\neq \emptyset$ and part (2) applies or $(b,c)\in J(R)^2$ or $(a,d)\in J(R)^2$. The case $(a,d)\in J(R)^2$ is entirely similar to the case $(b,c)\in J(R)^2$, hence we detail here the case when $(b,c)\in J(R)^2$. As $(b,c)\in J(R)^2$, $(a,d)\in Um(R^2)$; let $(e,f)\in R^2$ be such that $1=ae+df$. Then $\bigl(ae+bf+J(R),ce+df+J(R)\bigr)\in Um\bigl((R/J(R))^2\bigr)$, hence $(ae+bf,ce+df)\in Um(R^2)$; e.g., if $b=c=0$, then $\left[ 
\begin{array}{ccc}
a & 0 & f \\ 
0 & d & -e \\ 
-1 & 1 & 0\end{array}\right] $ is a simple extension of $A$.

\textbf{(4)} If $a,b,c,d,f,q\in R$ are such that $aq+df=1$, then a simple
extension of $\left[ 
\begin{array}{cc}
a & ab \\ 
ac & d\end{array}\right] $ is $\left[ 
\begin{array}{ccc}
a & ab & f \\ 
ac & d & -q+cf(1-b) \\ 
-1 & 1-b & 0\end{array}\right] $.
\end{example}

In three of these examples, the $(2,3)$ entry of the simple extensions is $-1 $, i.e., we can choose $e=1$. Such extensions relate to stable ranges 1 and 1.5 as follows.

\begin{corollary}\label{C10} 
The following properties hold:

\medskip \textbf{(1)} We have $sr(R)=1$ iff each upper triangular matrix $A\in Um\bigl(\mathbb{M}_2(R)\bigr)$ has a simple extension whose $(2,3)$
entry is $-1$.

\smallskip \textbf{(2)} We have $fsr(R)=1.5$ iff each upper 
triangular matrix $A\in Um\bigl(\mathbb{M}_2(R)\bigr)$ with nonzero $(1,1)$ entry has 
a simple extension whose $(2,3)$ entry is $-1$.

\smallskip \textbf{(3)} We have $asr(R)=1$ iff each upper 
triangular matrix $A\in Um\bigl(\mathbb{M}_2(R)\bigr)$ with $(1,1)$ entry not in $J(R)$ has  a simple extension whose $(2,3)$ entry is $-1$.
\end{corollary}

\begin{proof}
Let $A=\left[ 
\begin{array}{cc}
a & b \\ 
0 & d\end{array}\right] \in Um\bigl(\mathbb{M}_2(R)\bigr)$. If $a=0$, then $(b,d)\in Um(R^{2})$ and
from Equation (\ref{EQ1}) it follows that $A$ has a simple extension with
the $(2,3)$ entry $-1$ iff there exists $(f,t)\in R^{2}$ such that $bt+dft=1$ and hence iff there exists $f\in R$ such that $b+df\in
U(R)$. Thus all these matrices $A$ with $a=0$ have a simple extension with
the $(2,3)$ entry $-1$ iff $sr(R)=1$. 

Similarly, if $a\neq 0$ (resp.\ $a\notin J(R)$), then from Equation (\ref{EQ1}) it follows that $A$ has a simple extension with
the $(2,3)$ entry $-1$ iff there exists $(e,f,t)\in R^{3}$ such that $ae+bt+dft=1$ and hence iff there exists $f\in R$ such that $(b+df,a)\in
Um(R^2)$. From the definition of stable range 1.5 (resp.\ almost stable range 1) applied to $(b,d,a)\in Um(R^{3})$ it follows that all these matrices with $a\neq 0$ (resp.\ $a\notin J(R)$) have a simple extension with the $(2,3)$ entry $-1$ iff $fsr(R)=1.5$ (resp.\ $asr(R)=1$).
\end{proof}

\section{Proof of Theorem \protect\ref{TH1}}\label{S5}

\begin{proposition}
\label{PR5} Let $A\in Um\bigl(\mathbb{M}_2(R)\bigr)$. Then the following properties
hold:

\medskip \textbf{(1)} Assume that $\det (A)=0$. Then $A$ is simply extendable iff one of the three $R$-modules $\Im_A$, $\Ker_A$ and $\Coker_A$ is isomorphic to $R$ and iff $A$ is non-full.

\smallskip \textbf{(2)} The matrix $A$ is extendable iff its reduction modulo $R\det(A)$ is
non-full.
\end{proposition}

\begin{proof}
We prove part (1) using three (circular) implications; hence $\det(A)=0$. 

If $A$ is non-full, we write $A=\left[ 
\begin{array}{c}
l \\ 
m\end{array}\right] \left[ 
\begin{array}{cc}
o & q\end{array}\right] $. As $A\in Um\bigl(\mathbb{M}_2(R)\bigr)$, it follows that $(l,m),(o,q)\in
Um(R^{2})$; let $(e,f),(s,t)\in R^{2}$ be such that $el+fm=so+tq=1 $. Thus $\left[ 
\begin{array}{cc}
e & f\end{array}\right] A=\left[ 
\begin{array}{cc}
o & q\end{array}\right]$. Hence $A$ is simply
extendable by Theorem \ref{TH4}, a simple extension of it being $\left[ 
\begin{array}{ccc}
lo & lq & f \\ 
mo & mq & -e \\ 
-t & s & 0\end{array}\right] $.

If $A$ is simply extendable, it admits diagonal reduction by Theorem \ref{TH4}. Thus, $A$, being unimodular with $\det(A)=0$, is equivalent to $Diag(1,0)$, so $\Im(A)\cong R$. 

If one of the three $R$-modules $\Im_A$, $\Ker_A$ and $\Coker_A$ is isomorphic to $R$, then all of them are isomorphic to $R$ by Lemma \ref{L0}. As $\Im_A\cong R$, the $R$-linear $L_A:R^2\rightarrow R^2$ is a composite $R$-linear map $R^2\rightarrow R\rightarrow R^2$, hence $A$ is non-full. 

Thus part (1) holds. Part (2) follows from part (1) and Lemma \ref{L1}(1).
\end{proof}

We are now ready to prove Theorem \ref{TH1}. If $A\in Um\bigl(\mathbb M_2(R)\bigr)$ has zero determinant, then $\Ker_A$ and $\Im_A$ are projective $R$-modules of rank $1$ dual to each other and $\Ker_A\oplus\Im_A\cong R^2$ by Lemma \ref{L0}. Hence $(4)\Rightarrow (3)$. The equivalences $(1)\Leftrightarrow (2)\Leftrightarrow (3)$ follow from Proposition \ref{PR5}(1).

We show that $(2)\Rightarrow (4)$. Each projective $R$-module $P$ of rank $1$ generated by $2$ elements is isomorphic to $\Im_{A}$
for some idempotent $A\in \mathbb{M}_{2}(R) $ of rank $1$ and hence
unimodular of zero determinant. Assuming that (2) holds, we write $A=\left[\begin{array}{c}
l \\ 
m\end{array}\right] \left[ 
\begin{array}{cc}
o & q\end{array}\right] $ with $(l,m),(o,q)\in
Um(R^{2})$, hence $L_A$ is the composite of a surjective $R$-linear map $R^2\to R$ and an injective $R$-linear map $R\to R^2$. Thus $P\cong\Im_A\cong R$, hence $(2)\Rightarrow (4)$. Thus Theorem \ref{TH1} holds.

\begin{example}
\normalfont\label{EX7} For $A\in Um\bigl(\mathbb{M}_2(R)\bigr)$ with $\det
(A)=0$, we consider the $R$-submodule $K:=R\left[ 
\begin{array}{c}
-b \\ 
a\end{array}\right] +R\left[ 
\begin{array}{c}
-d \\ 
c\end{array}\right] $ of $\Ker_A$. For a maximal ideal $\mathfrak{m}$ of $R$, $K\not\subseteq\mathfrak{m}\mathbb{M}_{2}(R)$, hence $K\not\subseteq\mathfrak{m}\Ker_A$. Thus $\Ker_A=K$. For $(e,f)\in R^{2}$, we have $(ae+cf,be+df)\in Um(R^2)$ iff $\Ker_A$ is free having $e\left[ 
\begin{array}{c}
-b \\ 
a\end{array}\right] +f\left[ 
\begin{array}{c}
-d \\ 
c\end{array}\right] $ as a generator.
\end{example}

\begin{example}
\normalfont\label{EX8} Let $n\in\mathbb N$ and let $x_1,\ldots,x_n$ be indeterminates. Let $k\in \mathbb{N}$. Let $q:=4k+1$, $r:=2k+1$ and $\theta :=i\sqrt{q}\in \mathbb{C}$. We check
that $\mathbb{Z}[x_1,\ldots,x_n]$ \emph{is a }${\Pi}_{2}$\emph{\ ring which is not an} $E_{2}$ \emph{ring}. As $\mathbb{Z}[x_1,\ldots,x_n]$ is a UFD, it is also a Schreier domain and hence a ${\Pi}_{2}$
ring (see Section \ref{S1}). Thus it suffices to show that the matrix 
\begin{equation*}
A:=\left[ 
\begin{array}{cc}
r & 1-x_1 \\ 
1+x_1 & 2\end{array}\right] \in Um\bigl(\mathbb{M}_{2}(\mathbb{Z}[x_1,\ldots,x_n])\bigr)
\end{equation*}is not extendable. As $\det (A)=x_1^{2}+q$, based on Lemma \ref{L1}(1), it
suffices to show that the image $B:=\left[ 
\begin{array}{cc}
r & 1-\theta \\ 
1+\theta & 2\end{array}\right] \in Um\bigl(\mathbb{M}_{2}(\mathbb{Z}[\theta ])\bigr)$ of $A$, via the composite
homomorphism 
$\mathbb{Z}[x_1,\ldots,x_n]\rightarrow \mathbb{Z}[x_1,\ldots,x_n]/(x_1^{2}+q)\rightarrow \mathbb{Z}[\theta ]$
that maps $x_1$ to $\theta $ and $x_2, \ldots,x_n$ to $0$, is not simply extendable. An argument on norms shows that the element $2\in\mathbb{Z}[\theta ]$
is irreducible, i.e., $2=u (2u^{-1})$ with $u\in U(\mathbb Z[\theta ])$ are its only product decompositions. So, as $2u^{-1}$ divides
neither $1-\theta$ nor $1+\theta$, $B$ is full. So $B$ is not simply extendable by Proposition \ref{PR5}(1) and the integral domain $\mathbb{Z}[\theta ]$ is not a ${\Pi}_{2}$ ring. If $4k+1$ is square free, then $\mathbb{Z}[\theta ] $
is a Dedekind domain but not a \textsl{PID}.
\end{example}

\begin{remark}\label{rem2}
\normalfont
Statements (1) to (4) of Theorem \ref{TH4} are stable
under similarity (inner automorphisms of the $R$-algebra $\mathbb{M}_2(R)$)
but in general they are not stable under all $R$-algebra automorphisms
of $\mathbb{M}_2(R)$. To check this, let $R$ be such that there
exists an $R$-module $P$ such that $P\oplus P=R^2$ but $P\ncong R$; so $[P]\in {Pic}_{2}(R)[2]$. The idempotent $A$ of $\mathbb{M}_2(R)$ which is a projection of $R^2$ on the first
copy of $P$ along the second copy of $P$ satisfies $\Im_A=P$ and $\det(A)=0$, so it is not
simply extendable by Theorem \ref{TH1} but its image under the $R$-algebra automorphism $\mathbb{M}_2(R)={End}_R(P\oplus P)\cong \mathbb{M}_2(R)$ defined by $End_R(P)\cong R$ maps $A$ to $Diag(1,0)$. E.g., if $R$ is a Dedekind domain with ${Pic}(R)\cong\mathbb{Z}/2\mathbb{Z}$ (such as $\mathbb{Z}[\sqrt{-5}]$), then we can take $P$ to be a maximal ideal of $R$ with nontrivial class in ${Pic}(R)$.
\end{remark}

\section{Proofs of Theorems \ref{TH2} and \ref{TH3} and Corollary \ref{C2}}\label{S6}

We show that if $sr(R)\le 2$, then each extendable $A=\left[ 
\begin{array}{cc}
a & b \\ 
c & d\end{array}\right]\in Um\bigl(\mathbb M_2(R)\bigr)$ is simply extendable. 
Let $(e^{\prime},f^{\prime})\in R^2$ be such that 
$(ae^{\prime}+cf^{\prime},be^{\prime}+df^{\prime},ad-bc)\in Um(R^3)$ by Corollary \ref{C8}(1). 
Thus $(e^{\prime},f^{\prime},ad-bc)\in Um(R^3)$. As $sr(R)\le 2$, there exists $(r_1,r_2)\in R^2$ such that 
$$(e,f):=\bigl(e^{\prime}+(ad-bc)r_1,f^{\prime}+(ad-bc)r_2\bigr)\in
Um(R^2).$$ 
As the ideals of $R$ generated by $ae+cf,be+df,ad-bc$ and by $ae^{\prime}+cf^{\prime},be^{\prime}+df^{\prime},ad-bc$ are equal, 
we have $(ae+cf,be+df,ad-bc)\in Um(R^3)$. Hence $A$ is simply extendable by Corollary \ref{C8}(2). Thus Theorem \ref{TH2} holds.

\begin{example}
\normalfont\label{EX9} Let 
$R$ be an integral domain such that $sr(R)=3$ and each 
projective $R$-module $P$ with $P\oplus R\cong R^3$ is free. E.g., 
if $\kappa$ is a subfield of $\mathbb{R}$, then $sr(\kappa[x_{1},x_{2}])=3$ by \cite{vas}, Thm.\ 8 and Seshadri proved that all 
finitely generated projective modules over it are free (see \cite{ses}, Thm.; 
see also \cite{lan}, Ch.\ XXI, Sect.\ 3, Thm.\ 3.5 for Quillen--Suslin Theorem). Let 
$(a_{1},a_{2},b)\in Um(R^{3})$ be not reducible; thus $b\neq
0$. We have projective resolutions $0\rightarrow Rb\rightarrow R\rightarrow
R/Rb\rightarrow 0$ and $0\rightarrow P\xrightarrow{g}R^{2}\xrightarrow{f} R/Rb\rightarrow 0$, where the $R$-linear
map $f$ maps the elements of the standard basis of $R^{2}$ to $a_{1}+Rb$ and $a_{2}+Rb$, $P=\Ker(f)$ and $g$ is the inclusion. As $P$ is of the type mentioned (see Example \ref{EX4.1}), we identify $P=R^{2}$. Let $A\in \mathbb{M}_{2}(R)$ be such that $L_A=g:R^{2}=P\rightarrow R^{2}$; 
we have $A\in Um\bigl(\mathbb M_2(R)\bigr)$ and $R\det(A)=Rb$. Let $\bar A$ be the reduction of $A$ modulo 
$R\det(A)$. The $R/R\det(A)$-module $\Coker_A$ is isomorphic to 
$R/R\det(A)$ and to $\Coker_{\bar A}$. Thus $\bar A$ is simply extendable by Proposition \ref{PR5}(1). Hence $A$ is extendable by Lemma \ref{L1}(1). But $A$ is not simply extendable: if it were, then it would be equivalent to ${Diag}\bigl(1,\det(A)\bigr)$ (see Theorem \ref{TH4}) and it would follow from Example \ref{EX4.1} that $(a_{1},a_{2},b)\in Um(R^{3})$ is reducible, a contradiction.
\end{example}

To prove Theorem \ref{TH3}, let $R$ be an integral domain of dimension 1. To prove part (1), let $A\in Um\bigl(\mathbb{M}_2(R)\bigr)$ with $\det (A)\neq 0$. 
The ring $\bar{R}:=R/\det(A)R$ has dimension $0$ and hence $Pic(\bar{R})$ is trivial. 
From this and Theorem \ref{TH1} it follows that $\bar{R}$ is a $\Pi_2$ ring and therefore the reduction $\bar{A}$ of $A$ modulo $R\det (A)$ is simply extendable. From Theorem \ref{TH1} it follows that $\Im_{\bar{A}}\cong R/R\det (A)$. As $R$ is an integral domain, we have two projective resolutions $0\rightarrow R\det(A)\rightarrow R\rightarrow \Im_{\bar{A}}\rightarrow 0$ and $0\rightarrow \det(A)R^2\rightarrow\Im_A\rightarrow\Im_{\bar{A}}\rightarrow 0$ of $R$-modules. From this and Example \ref{EX4.1} applied to $b=\det(A)$ (recall that $sr(R)\le 2$), it follows that $A$ is equivalent to ${Diag}\bigl(1,\det (A)\bigr) $ and hence is simply extendable (see Theorem \ref{TH4}).

Part (2) holds as triangular matrices in $Um\bigl(\mathbb{M}_2(R)\bigr)$ have either two $0$ entries or nonzero determinants and thus are simply extendable by Example \ref{EX6}(3) or part (1). The first `iff' of part (3) follows from part (1) and definitions. The isomorphism classes of projective $R$-modules of rank 1 are the isomorphisms classes of nonzero ideals of $R$ which locally in the Zariski topology are principal; as for each $a\in R\setminus\{0\}$, $\dim(R/Ra)=0$ and hence $Pic(R/Ra)$ is trivial, all such nonzero ideals are generated by 2 elements. From this and Theorem \ref{TH1} it follows that $R$ is a $\Pi_2$ domain iff $Pic(R)$ is trivial. Hence part (3) holds. As \textsl{PID}s are precisely Dedekind domains with trivial Picard groups, part (4) follows from the second `iff' of part (3). Thus Theorem \ref{TH3} holds.

\begin{remark}\label{rem3}
\normalfont
The existence of an extension of a matrix in $Um\bigl(\mathbb{M}_2(R)\bigr)$ does not depend only on \emph{the set} of its entries. This is so
as, referring to Example \ref{EX8}, the matrix $\left[ 
\begin{array}{cc}
1+\theta & 1-\theta \\ 
r & 2\end{array}\right] \in Um\bigl(\mathbb{M}_{2}(\mathbb{Z}[\theta ])\bigr)$ has the same entries as $B$, has nonzero determinant, and \emph{it is} extendable (see Theorem \ref{TH3}(1)).
\end{remark}

To prove Corollary \ref{C2}, we assume that $asr(R)=1$. Each matrix $\left[ 
\begin{array}{cc}
a & b \\ 
0 & c\end{array}\right] \in Um\bigl(\mathbb{M}_2(R)\bigr)$ is simply extendable by Example \ref{EX6}(3) if $a\in J(R)$ and by Corollary \ref{C10}(3) if $a\notin J(R)$. So part (1) holds. To prove part (2), as $R$ is also a Hermite ring, each $A\in Um\bigl(\mathbb{M}_2(R)\bigr)$ is equivalent to a triangular matrix and hence from part (1) and Lemma \ref{L1}(3) it follows that $A$ is simply extendable and thus admits diagonal reduction by Theorem \ref{TH4}. As $\mathbb M_2(R)=RUm\bigl(\mathbb M_2(R)\bigr)$ it follows that each matrix in $\mathbb M_2(R)$ admits diagonal reduction, hence it is equivalent to a diagonal matrix. Thus $R$ is an EDR. Hence Corollary \ref{C2} holds.

\section{Explicit computations for integral domains}\label{S7}

Let $A:=\left[ 
\begin{array}{cc}
a & b \\ 
c & d\end{array}\right] \in Um\bigl(\mathbb{M}_2(R)\bigr)$ be such that we can write $a=ga^{\prime }$, $c=gc^{\prime }$, $b=hb^{\prime }$, $d=hd^{\prime }$ with $a^{\prime
},b^{\prime },c^{\prime },d^{\prime },g,h\in R$ and $(a^{\prime
},c^{\prime }),(b^{\prime },d^{\prime })\in Um(R^{2})$. We have $(g,h)\in
Um(R^{2})$. Let $e^{\prime },f^{\prime }\in R$ be such that $a^{\prime
}e^{\prime }+c^{\prime }f^{\prime }=1$. Let $l:=b^{\prime }c^{\prime
}-a^{\prime }d^{\prime }\in R$ and $m:=b^{\prime }e^{\prime }+d^{\prime
}f^{\prime }\in R $; note that $\det (A)=-ghl$. As $[\begin{array}{cc}
l & m\end{array}]=[\begin{array}{cc}
b^{\prime } & d^{\prime }\end{array}]\left[ 
\begin{array}{cc}
c^{\prime } & e^{\prime } \\ 
-a^{\prime } & f^{\prime }\end{array}\right] $, the matrix $\left[ 
\begin{array}{cc}
c^{\prime } & e^{\prime } \\ 
-a^{\prime } & f^{\prime }\end{array}\right] $ has determinant 1, and, due to 
$({b}^{\prime }{,d}^{\prime })\in Um(R^2)$, it follows that $(l,m)\in Um(R^2)$. 

Let $(e,f,w)\in R^3$ be such that $ae+cf=gw$. If $R$ is an
integral domain and $g\neq 0$, there exists $v\in R$ such that $(e,f)=(we^{\prime
}+c^{\prime }v,wf^{\prime }-a^{\prime }v)$, hence
\begin{equation*}
be+df=h(b^{\prime }e+d^{\prime }f)=h(b^{\prime }we^{\prime }+b^{\prime
}c^{\prime }v+d^{\prime }wf^{\prime }-a^{\prime }d^{\prime }v)=hw(b^{\prime}e^{\prime }+d^{\prime }f^{\prime })+hv(b^{\prime }c^{\prime }-a^{\prime}d^{\prime })
\end{equation*}is equal to $h(wm+vl)$. Thus $(ae+cf,be+df)=\bigl(gw,h(wm+vl)\bigr)$.

\begin{proposition}
\label{PR6} Let $R$ be an integral domain. Let $A=\left[ 
\begin{array}{cc}
a & b \\ 
c & d\end{array}\right] \in Um\bigl(\mathbb{M}_2(R)\bigr)$ be such that the above notation $g,h,a^{\prime },b^{\prime },c^{\prime },d^{\prime }$ applies and let $(e^{\prime },f^{\prime },l,m)\in R^{4}$ be obtained as above. 
We assume $g\neq 0$.

\medskip \textbf{(1)} The matrix $A$ is simply extendable iff
there exists $(w,v)\in R^2$ such that $(g,wm+vl),(w,hvl)\in Um(R^{2})$, in which
case a simple extension of $A$ is 
\begin{equation*}
\left[ 
\begin{array}{ccc}
a & b & wf^{\prime }-a^{\prime }v \\ 
c & d & -we^{\prime }-c^{\prime }v \\ 
-t & s & 0\end{array}\right]
\end{equation*}where $s,t\in R$ are such that $gws+h(wm+vl)t=1$.

\smallskip \textbf{(2)} If the intersection $\{(g,l),(g,m),(h,l),(h,m)\}\cap
Um(R^{2})$ is nonempty (e.g., if $hlm=0$), then $A$ is simply extendable
and $w,v\in R$ are given by formulas.
\end{proposition}

\begin{proof}
There exists $(e,f)\in R^2$ such that $(ae+cf,be+df)\in Um(R^{2})$ iff
there exists $(w,v)\in R^2$ such that $\bigl(gw,h(wm+vl)\bigr)\in Um(R^{2})$ (see above).
We have $\bigl(gw,h(wm+vl)\bigr)\in Um(R^{2})$ iff $\bigl(g,h(wm+vl)\bigr),\bigl(w,h(wm+vl)\bigr)\in Um(R^{2})$. As $(g,h)\in Um(R^{2})$, we have 
$\bigl(g,h(wm+vl)\bigr)\in Um(R^{2})$ iff $(g,wm+vl)\in Um(R^{2})$;
moreover, $Rw+Rh(wm+vl)=Rw+Rhvl$. Thus $\bigl(gw,h(wm+vl)\bigr)\in Um(R^{2})$ 
iff $(g,wm+vl),(w,hvl)\in Um(R^{2})$. Based on the `iff'
statements of this paragraph and Theorem \ref{TH4}, it follows that part
(1) holds.

To check part (2), we first notice that if $hlm=0$, then the intersection is nonempty; e.g., if $h=0$, then $g\in U(R)$
and hence $(g,l),(g,m)\in Um(R^{2})$. Based on part (1) it suffices to show
that in all four possible cases, we can choose $(w,v)\in R^2$ such that $(g,wm+vl),(w,hvl)\in Um(R^{2})$.

If $(g,l)\in Um(R^{2})$, for $(w,v):=(g,1)$ we have $Rg+R(wm+vl)=Rg+Rl=R$.
As $(g,h),(g,l)\in Um(R^{2})$, also $(w,hvl)=(g,hl)\in Um(R^{2})$.

If $(g,m)\in Um(R^{2})$, for $(w,v):=(1,0)$ we have $(g,wm+vl)=(g,m)\in
Um(R^{2})$ and $(w,hvl)=(1,0)\in Um(R^{2})$.

If $(h,m)\in Um(R^{2})$, then $(hl,m)\in Um(R^{2})$ and there exists $(w,v^{\prime })\in R^2$ such that $wm+hv^{\prime }l=1$; so 
$Rw+Rhv^{\prime}l=R $. For $v:=hv^{\prime }$ we have $wm+vl=1$ and 
so $(g,wm+vl)\in
Um(R^{2})$ and $(w,hvl)=(w,h^{2}v^{\prime }l)\in Um(R^{2})$ as $(w,hv^{\prime }l)\in Um(R^{2})$.

If $(h,l)\in Um(R^{2})$, let $(p,q)\in R^2$ be such that $1=lp+mq$. 
For $w:=hq+l$
and $v:=hp-m$ we compute $wm+vl=h(lp+mq)+ml-ml=h$, so $(g,wm+vl)=(g,h)\in Um(R^{2})$. Also, $Rw+Rh=R(hq+l)+Rh=Rl+Rh=R$ 
and $Rw+Rvl$ contains $wm+vl=h$ and hence it contains $Rw+Rh=R$. 
Thus $(w,vl)\in
Um(R^{2}) $. As $(w,h),(w,vl)\in Um(R^{2})$ it follows that $(w,hvl)\in
Um(R^{2})$.
\end{proof}

\begin{remark}\label{rem5}
\normalfont
We include a second proof of Corollary \ref{C2}(2) for B\'{e}zout domains. Assume $R$ is a B\'{e}zout domain with $asr(R)=1$. Based on the equivalence of statements (1) and (2) of Theorem \ref{TH4}, it suffices to
show that each matrix $A=\left[ 
\begin{array}{cc}
a & b \\ 
c & d\end{array}\right] \in Um\bigl(\mathbb{M}_2(R)\bigr)$ is simply extendable. As $R$ is a 
Hermite domain, the notation of this section applies. As $(g,h)\in Um(R^2)$, 
by the symmetry between the pairs $(a,c)$ and $(b,d)$, 
we can assume that $g\notin J(R)$. As $(m,l,g)\in
J_{3}(R)$ and $asr(R)=1$, there exists $v\in R$ such that $(m+lv,g)\in Um(R^2)$.
We take $w:=1$. Hence $(w,hvl)\in Um(R^2)$ and $(g,wm+lv)=(g,m+lv)\in Um(R^2)$.
From Proposition \ref{PR6}(1) it follows that $A $ is simply extendable.
\end{remark}

\begin{example}\normalfont\label{EX10}
For $A=\left[ 
\begin{array}{cc}
a & b \\ 
c & d\end{array}\right]\in Um\bigl(\mathbb M_2(\mathbb Z)\bigr)$, its simple extensions $A^+_{(e,f,s,t)}=\left[ 
\begin{array}{ccc}
a & b & f \\ 
c & d & -e \\ 
-t & s & 0\end{array}\right]$ are parameterized by the set (see Equation (\ref{EQ1}))
$$\gamma_A:=\{(e,f,s,t)\in \mathbb{Z}^{4}|a(es)+b(et)+c(fs)+d(ft)=1\}.$$
The below examples were (initially) exemplified using a code written for $R=\mathbb Z$ by the second author. We have $\nu_{A^+_{(e,f,s,t)}}=\chi_{A^+_{(e,f,s,t)}}^{\prime}(0)=\det(A)+es+ft$.

(1) If $a=0$ and $d=1+b+c$, then we can take $A^{+}=\left[ 
\begin{array}{ccc}
0 & b & -1 \\ 
c & 1+b+c & -1 \\ 
1 & 1 & 0\end{array}\right].$

Concretely, suppose $(b,c)=(3,2)$; then $d=6$, $\det(A)=-6$, $\gamma_A=\{(e,f,s,t)\in \mathbb{Z}^{4}|3et+2fs+6ft=1\}$ and $\nu(A)=\{-6+es+ft|(e,f,s,t)\in\gamma_A\}$.

\smallskip To solve the equation $3et+2fs+6ft=1$, let $w:=et+2ft$. We get 
$2fs+3w=1$ with general solution $fs=-1+3k$, $w=1-2k$, where 
$k\in\mathbb{Z}$. The general solution of the equation $et+2ft=1-2k$ is 
$et=2k-1-2l$ and $ft=1-2k+l$, where $l\in\mathbb{Z}$. Let $m:=l-2k+1$. 
It follows that $ft=m,\; et=1-2k-2m,\; fs=-1+3k$
and the only constraint is that $ft=m$ divides $etfs=(3k-1)(1-2k-2m)$, i.e., divides $(3k-1)(2k-1)$. As $es+ft=m+\dfrac{-6k^2+5k-1}{m}-6k+2$, it follows that
\begin{equation*}
\nu(A)=\{-4+m-6k+\dfrac{-6k^2+5k-1}{m}|(m,k)\in\mathbb{Z}^2,\; m\;{\textup{divides}}\; -6k^2+5k-1\}.
\end{equation*}


(2) If $c=0$ and $d=1-a+b$, then we can take $A^{+}=\left[ 
\begin{array}{ccc}
a & b & -1 \\ 
0 & 1-a+b & -1 \\ 
1 & 1 & 0\end{array}\right] $ (cf.\ Corollary \ref{C9}(2): for simple extensions of upper triangular matrices with nonzero $(1,1)$ entries over rings $R$ with $fsr(R)=1.5$, we can choose $e=1$).

Concretely, suppose $(a,b)=(6,-10)$; hence $d=15$ and $\det(A)=-90$. Thus $\gamma_A=\{(e,f,s,t)\in \mathbb{Z}^{4}|6es-10et-15ft=1\}$.

To solve the equation $6es-10et-15ft=1$, let $w:=2et+3ft$. We get $6es-5w
=1$ with general solution $es=1+5k$, $w =1+6k$, where $k\in\mathbb{Z}$. Then 
$2et+3ft=1+6k$ has the general solution $et=-1-6k+3l$, $ft=1+6k-2l$, where $l\in\mathbb{Z}$. Let $m:=l-2k$. It follows that 
$es=1+5k,\; et=-1+3m,\; ft=1+2k-2m$ are subject to the only constraint that $et=-1+3m$ divides $esft=(1+5k)(1+2k-2m)$, i.e., it divides $(1+5k)(2k+m)$. As $es+ft=2+7k-2m$, it follows that 
\begin{equation*}
\nu(A)=\{-88+7k-2m|(m,k)\in\mathbb{Z}^2,\; -1+3m\;{\textup{divides}}\;
(1+5k)(2+2k+m)\}.
\end{equation*}
For $m=0$ (resp.\ $m=1$) it follows that $\nu(A)\supseteq 3+7\mathbb Z$ (resp.\ $\nu(A)\supseteq 1+14\mathbb{Z}$).

The matrices $Diag(6,-15),A_0:=A\in\mathbb{M}_2(\mathbb{Z}[\dfrac{1}{21}])$ are similar and one checks that $\nu(A_0)=-90+\dfrac{1}{6}+\dfrac{7}{2}+2\mathbb{Z}[\dfrac{1}{21}]=-86-\dfrac{1}{3}+2\mathbb{Z}[\dfrac{1}{21}]$.

(3) If $A=$ $\left[ 
\begin{array}{cc}
15 & 6 \\ 
10 & 14\end{array}\right] $, we can take $A^{+}=\left[ 
\begin{array}{ccc}
15 & 6 & -2 \\ 
10 & 14 & 1 \\ 
-1 & -1 & 0\end{array}\right] $; indeed we have $\det(A^+)=1\cdot 15-1\cdot 6+2\cdot 10-2\cdot 14=1$. The
entries of $A$ use double products of the primes $2,3,5,7$. We have 
$\gamma_A=\{(e,f,s,t)\in \mathbb{Z}^{4}|15es+6et+10fs+14ft=1\}$, $\det(A)=150$, $\nu_{A^{+}}=149$, and $\nu(A)=\{150+es+ft|(e,f,s,t)\in\gamma_A\}$.

To solve the equation $15es+6et+10fs+14ft=1$, let $x:=5es+2et$ and $y:=5fs+7ft$, so we get $3x+2y =1$ with general solution $x=1+2k$, $y =-1-3k$, where $k\in\mathbb{Z}$. Then $5es+2et=1+2k$ has the general solution $es=1+2k+2l$, $et=-2(1+2k)-5l$, where $l\in\mathbb{Z}$, and $5fs+7ft=-1-3k$
has the general solution $fs=3(-1-3k)+7r$, $ft=-2(-1-3k)-5r$, where $r\in\mathbb{Z}$. Let $o:=k+l$, so $et=-2+4l-4o-5l=:q$. Thus $es=1+2o$, $et=q$, 
$l=-2-q-4o$, $k=2+q+5o$, therefore $fs=-3-18-9q-45o+7r$ and 
$ft=2+12+6q+30o-5r$. Let $p:=r-6o-q-3$. Thus $ft=-5p+q-1$ and 
$fs=-3o-2q+7p$. As $(es)(ft)=(et)(fs)$, we have an identity $
(1+2o)(q-1-5p)=q(7p-3o-2q)$,
which can be rewritten as $o(-2-10p+5q)=1+5p-q-2q^2+7pq$
and which for $q=2p$ becomes $-2o=6p^2+3p+1$, requiring $p$ to be odd. For $q=2p$, $es+ft=2o+q-5p$ becomes $2o-3p=-6p^2-6p-1$. Thus
\begin{equation*}
\nu(A)=\{150+2o+q-5p|(o,q,p)\in\mathbb{Z}^3,\;o(-2-10p+5q)=1+5p-q-2q^2+7pq\}
\end{equation*}
contains the set $\{150-6p^2-6p-1|p-1\in 2\mathbb{Z}\}$.

(4) If $A=\left[ 
\begin{array}{cc}
30 & 42 \\ 
70 & 105\end{array}\right] $, we can take $A^{+}=\left[ 
\begin{array}{ccc}
30 & 42 & 1 \\ 
70 & 105 & 3 \\ 
1 & 1 & 0\end{array}\right] $; indeed we have $\det(A^+)=-3\cdot 30+3\cdot 42+1\cdot 70-1\cdot 105=1$. The
entries of $A$ use triple products of the primes $2,3,5,7$, hence they are not pairwise coprime.
\end{example}

\bigskip \noindent \textbf{Acknowledgement.} The third author would like to
thank SUNY Binghamton for good working conditions and Ofer Gabber for pointing out a glitch in one of the references. 
Competing interests: the authors declare none.

\hbox{} \hbox{Grigore C\u{a}lug\u{a}reanu\;\;\;E-mail: calu@math.ubbcluj.ro}
\hbox{Address: Department of Mathematics, Babe\c{s}-Bolyai
University,} 
\hbox{1 Mihail Kogălniceanu Street, Cluj-Napoca 400084, Romania.}

\hbox{} \hbox{Horia F.\ Pop\;\;\;E-mail: horia.pop@ubbcluj.ro} 
\hbox{Address:
Department of Computer Science, Babe\c{s}-Bolyai
University,} 
\hbox{1 Mihail Kogălniceanu Street, Cluj-Napoca 400084, Romania.}

\hbox{} \hbox{Adrian Vasiu,\;\;\;E-mail: avasiu@binghamton.edu} 
\hbox{Address:
Department of Mathematics and Statistics, Binghamton University,} 
\hbox{P.\ O.\ Box
6000, Binghamton, New York 13902-6000, U.S.A.}


\begin{thebibliography}{99}

\bibitem{bas} H.\ Bass \textsl{K-theory and stable algebra.} Inst.\ Hautes 
\'{E}tudes Sci.\ Publ.\ Math.\ \textbf{22} (1964), 5--60.

\bibitem{coh} P.\ M.\ Cohn \textsl{B\'{e}zout rings and their subrings.} Math.\ 
Proc.\ Camb.\ Philos.\ Soc.\ \textbf{64} (1968), 251--264.

\bibitem{CP} G.\ C\u{a}lug\u{a}reanu, H.\ F.\ Pop \textsl{On zero determinant
matrices that are full.} Math.\ Panon.\ New Series \textbf{27} /NS 
\textbf{1}/ (2021), no.\ 2, 81--88.

\bibitem{FL} H.\ K.\ Farahat, W.\ Ledermann \textsl{Matrices with prescribed
characteristic polynomials.} Proc.\ Edinb.\ Math.\ Soc.\ \textbf{11}
1958/1959, 143--146.

\bibitem{GH1} L.\ Gillman, M.\ Henriksen \textsl{Rings of continuous functions
in which every finitely generated ideal is principal.} Trans.\ Amer.\ Math.
Soc.\ \textbf{82} (1956), 366--391.

\bibitem{hei} R.\ C.\ Heitmann \textsl{Generating non-Noetherian modules efficiently.}
Michigan Math.\ J.\ \textbf{31} (1984), no.\ 2, 167--180.

\bibitem{her} D.\ Hershkowitz \textsl{Existence of matrices with prescribed
eigenvalues and entries.} Linear Multilinear Algebra \textbf{14} (1983),
no.\ 4, 315--342.

\bibitem{kap} I.\ Kaplansky \textsl{Elementary divisors and modules.} Trans.
Amer.\ Math.\ Soc.\ \textbf{66} (1949), 464--491.

\bibitem{lan} S.\ Lang \textsl{Algebra. Revised third edition.} Grad.\ Texts
in Math., \textbf{211}, Springer-Verlag, New York, 2002.

\bibitem{LLS} M.\ D.\ Larsen, W.\ J.\ Lewis, T.\ S.\ Shores \textsl{Elementary
divisor rings and finitely presented modules.} Trans.\ Amer.\ Math.\ Soc.\ 
\textbf{187} (1974), 231--248.

\bibitem{LL} J.\ Liu, D.\ Li \textsl{The generalized Serre problem over
K-Hermite rings.} J.\ Syst.\ Sci.\ Complex.\ \textbf{30} (2017), no.\ 2, 510--518.

\bibitem{lor} D.\ Lorenzini \textsl{Elementary divisor domains and B\'{e}zout
domains.} J.\ Algebra \textbf{371} (2012), 609--619.

\bibitem{mcg1} W.\ W.\ McGovern \textsl{B\'{e}zout rings with almost stable
range 1.} J.\ Pure Appl.\ Algebra \textbf{212} (2008), no.\ 2, 340--348.

\bibitem{MR} S.\ McAdam, D.\ E.\ Rush \textsl{Schreier rings.} Bull.\ Lond.
Math.\ Soc.\ \textbf{10}, 1 (1978), 77--80.

\bibitem{MM} P.\ Menal, J.\ Moncasi \textsl{On regular rings with stable range
2.} J.\ Pure Appl.\ Algebra \textbf{24} (1982), no.\ 1, 25--40.

\bibitem{ses} C.\ S.\ Seshadri \textsl{Triviality of vector bundles over the affine 
space $K^2$.} Proc.\ Nat.\ Acad.\ Sci.\ U.S.A.\ \textbf{44} (1958), 456--458.

\bibitem{shc1} V.\ P.\ Shchedryk \textsl{Some properties of primitive matrices
over B\'{e}zout B-domain.} Algebra Discrete Math.\ 2005, no.\ 2, 46--57.

\bibitem{vas} L.\ N.\ Vaser\v{s}te\u{\i}n \textsl{The stable range of rings
and the dimension of topological spaces.} (Russian) Funktsional.\ Anal.\ i Prilo\v{z}en.\ \textbf{5} (1971), no.\ 2, 17--27. [English Translation: Funct.
Anal.\ Appl.\ \textbf{5} (1971), 102--110.]

\bibitem{VS} L.\ N.\ Vaser\v{s}te\u{\i}n, A.\ A.\ Suslin  \textsl{Serre's problem 
on projective modules over polynomial rings, and algebraic K-theory.} (Russian) 
Izv.\ Akad.\ Nauk SSSR Ser.\ Mat.\ \textbf{40} (1976), no.\ 5, 993--1054, 1199.

\bibitem{WW} R.\ Wiegand, S.\ Wiegand \textsl{Finitely generated modules over B\'{e}zout rings.} Pacific J.\ Math.\ \textbf{58} (1975), no.\ 2, 655--664.

\bibitem{zab1} B.\ V.\ Zabavsky \textsl{Diagonal reduction of matrices over
rings.} Math.\ Stud.\ Monog.\ Ser., Vol.\ \textbf{16}, VNTL Publishers, L'viv, 2012.
251 pp.

\bibitem{zaf} M.\ Zafrullah \textsl{On a property of pre-Schreier domains.}
Comm.\ Algebra \textbf{15} (1987), no.\ 9, 1895--1920.
\end{thebibliography}
\end{document}